\documentclass[12pt]{article}
\usepackage{latexsym,amssymb,amsmath}
\usepackage{amsthm, amstext}
\usepackage{array, amsfonts, mathrsfs}
\usepackage{mathrsfs}
\usepackage[dvips]{graphicx,color}
\usepackage{epsf}

\usepackage {amsfonts,amsmath}
\allowdisplaybreaks

\usepackage [dvips] {graphicx}
\setkeys {Gin} {draft=true} \setkeys {Gin} {draft=false}

\renewcommand {\d} {\,{\rm d}}    
\newcommand {\ddd} {{\rm d}}      
\newcommand {\rsp} {\mathbb{R}}   
\newcommand {\NB}  {\nobreak}
\newcommand {\supp}{\mathop{\rm supp}\nolimits}
\newcommand {\EQF} [1] {\mathbin {\,{\stackrel {\mbox {\footnotesize #1}} {=}}\,}}
\newcommand {\Widetilde} [1] {{\,\,\,\widetilde{\!\!\!#1\!\!\!} \,\,\,}}
\newcommand {\ZZ} {\mkern 2mu}
\newcommand {\la} {\langle}
\newcommand {\ra} {\rangle}

\begin{document}

\author {M.I.Belishev%
\thanks{The original version is published in
\newblock {\em Ill-Posed and Inverse Problems.
S.I.Kabanikhin and V.G.Romanov (Eds). VSP}, 2002, 55--72.}}

\title {How to see waves under the Earth surface \\ (the BC-method for geophysicists)}

\date {}

\maketitle

\begin{abstract}
The BC-method is an approach to inverse problems based on their
relations to the boundary control theory. The paper provides a
simple and physically transparent description of the method in the
case of dynamical inverse data given at a portion of the boundary.
\end{abstract}

\section*{Introduction}

The BC-method is an approach to inverse problems based on their
relations to the boundary control theory. The goal of this paper
is to provide a maximally simple and physically transparent
description of the approach in the case of dynamical inverse data
given at a portion of the boundary. It is the variant of the
BC-method, which is most promising for applications in geophysics
and acoustics.

A specific feature of the BC-method is that it recovers
(visualizes) the waves, whereas recovering the parameters, in a
sense, turns out to be a ``free addition''. The method is
time-optimal: the longer is the time of measurements at the
surface, the greater is the depth of recovering the waves and
parameters.

\section{Geometry of rays}
\subsection{c-metric}

We denote by $x= (x^1,x^2,x^3)$ the points of the space~$\rsp^3$.

Let $\Omega \subset \rsp^3$ be a domain (possibly, unbounded) with
a smooth boundary~$\Gamma$, and let $c=c(x)$ be the {\it speed\/}
of waves propagating in~$\Omega$. We assume that $0< c_*\leqslant
c(x) \leqslant c^*<\infty$ holds for some constants $c_*,\,c^*$.

The speed determines the c-{\it metric\/} with a distance
\begin{equation}\label{1.1}
\tau (x,y):= \inf \int_x^y \frac {|\d s|} {c(s)} \,,
\end{equation}
where the infimum is taken over the set of smooth curves
connecting $x$ with~$y$; thus $\tau (x,y)$ is the travel time
needed for a wave initiated at $x$ to reach~$y$. The geodesics of
the c-metric are the curves realizing the infimum in (\ref{1.1}).
If $c(x)\equiv\NB 1$, the c-metric is Euclidean: $\tau
(x,y)=|x-\NB y|$ and the corresponding geodesics are straight
lines.

Let $\sigma \subset \Gamma$ be a portion (open subset) of the
boundary. Introduce an {\it eikonal}
$$
\tau (x):=\min\limits_{y\in \sigma} \tau (x,y)
$$
and its level surfaces ({\it fronts\/})
$$
\Gamma^\xi :=\{ x\in \Omega \mid \tau (x)=\xi \} ,\quad \xi \geq 0.
$$
In dynamics, $\Gamma^\xi$ is a forward front (at the moment $t=\NB
\xi$) of the wave  initiated at~$\sigma$ and moving into~$\Omega$,
and $\tau (x)$ is the time needed for the wave to reach the
point~$x$. Let
$$
\Omega^\xi :=\{ x\in \Omega \mid \tau (x)<\xi \}
$$
be the subdomain filled with waves at $t=\NB \xi$. This subdomain is bounded by
$\Gamma^\xi$ and~$\Gamma$.

\subsection{Ray coordinates}

Fix a point $\gamma \in \sigma$; let $r_\gamma$ be the geodesic of
the c-metric starting from $\gamma$ orthogonally to~$\sigma$. We
shall denote by $r_\gamma [0,\xi ]$ a segment of~$r_\gamma$ of the
c-length~$\xi$; let $x(\gamma ,\xi )$ be a second endpoint of this
segment.

Fix $T>0$; the rays starting from $\sigma$ cover the subdomain
$$
B^T:=\bigcup\limits_{\gamma \in \sigma} \bigcup\limits_{0\leq \xi \leq T}
x(\gamma ,\xi )
$$
which we call {\it a tube}, $\sigma$ being its bottom. The tube
can be represented as a collection of the cross-sections
$$
\sigma^\xi :=B^T \cap \Gamma^\xi =
\bigcup\limits_{\gamma \in \sigma} x(\gamma ,\xi ),
$$
each $\sigma^\xi$ being a part of the front $\Gamma^\xi$ lighted with
rays. Notice that $\sigma^0$ coincides with~$\sigma$.

As is known, if $T$ is not too large, then the families of
rays~$r_\gamma$ and fronts~$\sigma^\xi$ are regular. This enables
one to introduce a special coordinate system in~$B^T$. The pair
$(\gamma ,\xi )$ is said to be the {\it ray coordinates\/} of a
point $x\in\NB B^T$ if $x$ lies at the ray $r_\gamma$ and the
front $\sigma^\xi$ (i.e., if $x=x(\gamma ,\xi )$). Notice that
$x(\gamma,0)=\NB \gamma$.

\begin{figure}[htbp]
\centering
\includegraphics {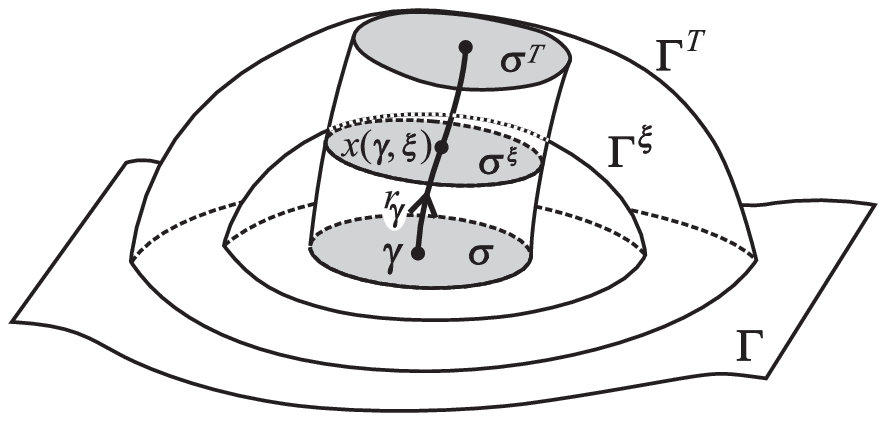}\\[1.5ex]
\small\bf Fig.~1. Ray coordinates
\end{figure}
Actually, since $\gamma$ is a point of the boundary, speaking
about coordinates we should provide $\gamma$ with local
coordinates $\gamma^1,\gamma^2$ on~$\sigma$ but such details are
not substantial for what follows.

\subsection{Ray divergence}

An important characteristic of the ray coordinates is a {\it ray
divergence}, which plays the role of an amplitude factor in the
Geometric Optics formulas.

Fix $\gamma \in \sigma$ and denote by $\sigma_\varepsilon (\gamma
)$ the intersection of $\sigma$ with the ball of small radius
$\varepsilon$ centered at~$\gamma$. Consider the tube
$B^T_\varepsilon$ formed by rays starting from~
$\sigma_\varepsilon (\gamma)$; let $\sigma^\xi_\varepsilon
(\gamma):=B^T_\varepsilon \cap\NB \Gamma^\xi$, and let $|\sigma|$
be an area of~$\sigma$. A function of the ray coordinates
$$
J(\gamma ,\xi ):=\lim\limits_{\varepsilon \to 0}
\frac {|\sigma^\xi_\varepsilon (\gamma)|} {|\sigma_\varepsilon (\gamma)|}
$$
is called the ray divergence at the point $x(\gamma ,\xi )$
(see, e.g.,~[\ref{BB}]).

\section{Propagation of jumps}
\subsection{Dynamical system with boundary control}

We denote by $Q^T:=\Omega \times\NB (0,T)$ a space-time cylinder;
let $\Sigma^T:=\Gamma \times\NB [0,T]$ be its lateral surface. The
part $\Sigma^T_\sigma :=\sigma \times\NB [0,T]$ of the lateral
surface is referred to as the {\it screen\/}; the meaning of this
term will be clarified later.

Consider the boundary initial-value problem
\begin{align}\label{2.1}
& u_{tt}-c^2\,[\Delta u-qu]=0 \qquad \mbox{in \,} Q^T, \\
\label{2.2} &
u|_{t=0}=u_t|_{t=0}=0 \qquad \mbox{in \,} \Omega ,\\
\label{2.3} & u=f \qquad \mbox{on \,} \Sigma^T,
\end{align}
where $q=q(x)$ is a {\it potential\/} in~$\Omega$ and $f=f(\gamma
,t)$ is a {\it boundary control}. Let $u=u^f(x,t)$ be a solution.
This solution describes a wave produced by the control~$f$ and
moving into the domain from the boundary.

In what follows we deal only with controls supported at the screen: $\supp
f\subset\NB \Sigma^T_\sigma$ (i.e., $f(\gamma ,t)=\NB 0$ for $\gamma \in
\Gamma \setminus\NB \sigma$ and for any~$t$). At the moment $t=\NB\xi$ the
forward front of the corresponding wave~$u^f$ coincides with the
surface~$\Gamma^\xi$ and the wave is supported in the subdomain~$\Omega^\xi$:
\begin{equation}\label{2.4}
\supp u^f(\cdot ,\xi )\subset \Omega^\xi,
\end{equation}
i.e., $u^f(x,\xi )= 0$ in $\Omega \setminus \Omega^\xi$. At the
final moment the wave captures the subdomain~$\Omega^T$.

If the control is switched on with a delay $T-\NB \xi$ \,(i.e.,
$f|_{0\leqslant t\leqslant T-\xi} =\NB 0$, so that $f$ acts
$\xi$~units of time), then at the final moment $t=\NB T$ the wave
$u^f(\cdot ,T)$ is  supported in~$\Omega^\xi$.

\subsection{Geometric Optics formula}

Take a smooth control $f=f(\gamma ,t)$ and fix an intermediate
moment~$t = T-\xi$,
$0<\xi <\NB T$. Let
$$
f_\xi (\gamma ,t):=\begin{cases}
0, & 0\leqslant t<T-\xi ,\\
f(\gamma ,t),& T-\xi \leqslant t\leqslant T,
\end{cases}
$$
be the truncated control. An important fact is that the truncation
violates the smoothness and leads to appearance of a jump at the moment
$t=T-\NB \xi$:
\begin{equation}\label{2.5}
f_\xi (\gamma ,T-\xi +0)=f(\gamma ,T-\xi ),
\end{equation}
so that the amplitude of the jump at the point
$(\gamma ,T-\NB \xi )\in\NB \Sigma^T_\sigma$
is equal to the value of the original control~$f$ at
this point.

The well-known fact is that a discontinuous control produces a discontinuous
wave with a jump propagating along rays with speed $c(x)$.
In our case the wave $u^{f_\xi}$ produced by the truncated control
enters into the domain with the initial jump
\begin{equation}\label{2.6}
u^{f_\xi}(x(\gamma ,0),T-\xi +0) \EQF {(\ref{2.3})}
f_\xi (\gamma ,T-\xi +0) \EQF {(\ref{2.5})} f(\gamma ,T-\xi ).
\end{equation}
Then this jump propagates along the ray $r_\gamma$ and at the final
moment $t=\NB T$ arrives at the point $x(\gamma ,\xi )$ (see Fig.~2).
Its final amplitude
$u^{f_\xi}(x(\gamma ,\xi -\NB 0),T)$ is connected with the initial one through the
well-known Geometric Optics law:
$$
\frac {u^{f_\xi}(x(\gamma ,\xi -0),T) } {u^{f_\xi}(x(\gamma ,0),T-\xi +0)}
=  \sqrt{\frac {c(x(\gamma ,\xi ))} {J(\gamma ,\xi )}} \Bigg/
 \sqrt{\frac {c(x(\gamma ,0 ))} {J(\gamma ,0 )}}
$$
(see, e.g.,~[\ref{BB}]), which implies
\begin{equation}\label{2.7}
{u^{f}(x(\gamma ,\xi -0),T) } \overset{(\ref{2.6})}=\sqrt { \frac
{c(x(\gamma ,\xi )) \, J(\gamma ,0 )}
 {c(x(\gamma ,0 )) \, J(\gamma ,\xi )} }
\,f(\gamma ,T-\xi ).
\end{equation}

\begin{figure}[htbp]
\centering
\includegraphics {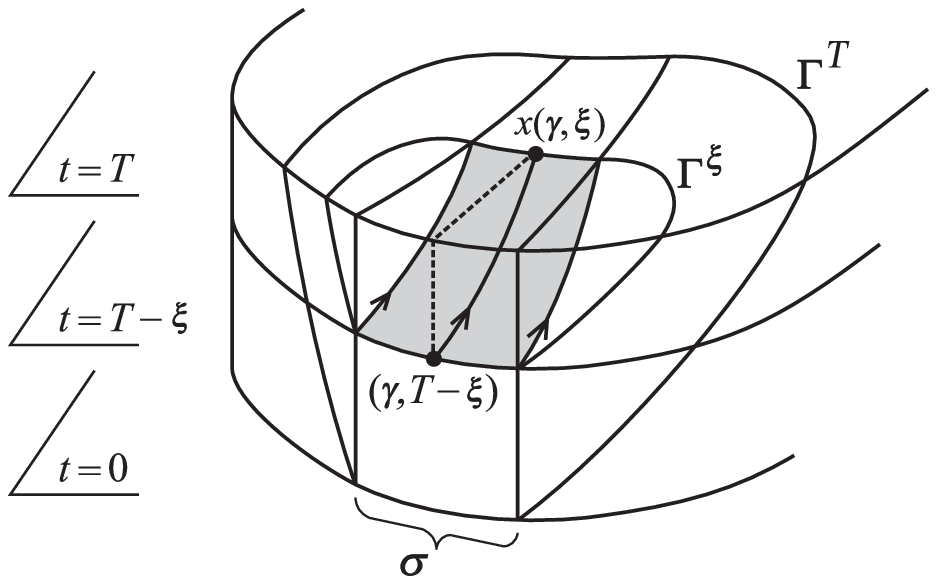}\\[1.5ex]
\small\bf Fig.~2. The jump of $u^{f_\xi}$
\end{figure}
Thus, at the final moment $t=\NB T$, the wave $u^{f_\xi}$ occupies
the subdomain~$\Omega^\xi$ and its forward front takes up the
position~$\Gamma^\xi$. At the part~$\sigma^\xi$ of the front,
which is lighted by the rays~$r_\gamma$, the wave has the jump,
whose amplitude can be calculated by~(\ref{2.7}).

Note in addition that the part $\Gamma^\xi \setminus\NB
\sigma^\xi$ of the front may contain singularities of rather
complicated structure. Fortunately, this part of the front plays
no role in further considerations.

\subsection{Dual system}

A dynamical system
\begin{align}
\label{2.8}& v_{tt}-c^2\,[\Delta v-qv]=0 \qquad \mbox{in \,} Q^T, \\
\label{2.9}& v|_{t=T}=0,\quad  v_t|_{t=T}=y \qquad \mbox{in \,} \Omega, \\
\label{2.10}&  v=0 \qquad \mbox{on \,} \Sigma^T
\end{align}
is called {\it dual\/} to system (\ref{2.1})--(\ref{2.3}). We
denote its solution by $v=v^y(x,t)$. Note that the Cauchy
data~(\ref{2.9}) are posed at the final moment; this is not
essential in view of the well-known time reversibility of
equation~(\ref{2.8}).

The solution $v^y$ describes a wave, which is produced by a speed
perturbation $y=y(x)$ and propagates in the domain with the
rigidly fixed boundary $\Gamma$. A function $\partial v^y (\gamma
,t) /
\partial \nu$ ($\nu =\nu (\gamma )$ is an outward normal
to~$\Gamma$) is proportional to the force, which appears as a
result of interaction between the wave~$v^y$ and the boundary at
the point~$\gamma$ at the moment~$t$.

Let us derive a relation between the solutions $u^f$ and $v^y$ of
the systems \mbox {(\ref{2.1})--(\ref{2.3})} and \mbox
{(\ref{2.8})--(\ref{2.10})}. It is the relation, which motivates
the term ``dual'': for any control~$f$ and perturbation~$y$ the
relation
\begin{equation}\label{2.11}
\int_\Omega \frac {\d x} {c^2(x)} \, u^f(x,T)\,y(x)=
\int_{\Sigma^T} \d\Gamma {\d t} \,\ZZ f(\gamma ,t)\, \frac
{\partial v^y} {\partial \nu} (\gamma ,t).
\end{equation}
holds. Indeed, integrating by parts in the equalities
\begin{align*}
0={}& \int_0^T \ddd t \int_\Omega \frac {\ddd x} {c^2(x)}
\, \big\{ u_{tt}^f (x,t) -c^2(x) \,\big[\Delta u^f(x,t)-q(x)\,u^f(x,t) \big] \big\}
\, v^y(x,t) \\
={}&\int_\Omega \frac {\ddd x} {c^2(x)}
\, \big[u_{t}^f (x,t)\,v^y(x,t) - u^f(x,t)\,v^y_t(x,t) \big] \big|_{t=0}^{t=T} \\
&{} -\int_0^T \ddd t \int_\Gamma \ddd\Gamma \,\Big[\frac
{\partial u^f} {\partial \nu} (\gamma ,t)\,v^y(\gamma ,t)-
u^f(\gamma ,t)\,\frac {\partial v^y} {\partial \nu} (\gamma ,t)\Big] \\
&{} +\int_0^T \ddd t \int_\Omega \frac {\ddd x} {c^2(x)}
\, u^f(x,t) \, \big\{ v_{tt}^y (x,t) -c^2(x)\, \big[\Delta v^y(x,t)-q(x)
\,v^y(x,t) \big] \big\} \\
\intertext {(see (2.2), (2.3), (2.9), and (2.10))}
={} & -\int_\Omega \frac {\ddd x} {c^2(x)}
\, u^f (x,T)\,y(x) +
\int_0^T \ddd t \int_\Gamma \ddd\Gamma
\,\ZZ f(\gamma ,t) \,\frac {\partial v^y} {\partial \nu} (\gamma ,t)
\end{align*}
we obtain (\ref{2.11}). To justify these calculations, one needs
to begin with taking the control~$f$ smooth and vanishing near
$t=\NB 0$, so that the wave $u^f(x,t)$ turns out to be smooth
in~$\Omega$ and vanishing near its forward front. Due to the
latter the surface integrals over the forward front vanish when we
integrate by parts. Then the final result (\ref{2.11}) is extended
to a wide class of (possibly discontinuous) controls~$f$. Notice
also that the integral in the right-hand side of (\ref{2.11}) is
in fact taken over $\Sigma^T_\sigma$ because the controls are
supported on the screen.

In what follows we write (\ref{2.11}) in a convenient symbolic
form. Introduce a {\it control operator\/} $W^T$ associated with
the system \mbox{(\ref{2.1})--(\ref{2.3})}, which maps controls to
a waves:
$$
W^T f:=u^f(\cdot ,T) \qquad \mbox{in \,} \Omega.
$$
It can be considered as an operator, which creates the waves.
Introduce also the {\it observation operator\/} $O^T$ associated
with the dual system, which maps perturbation $y$ to the force
produced by $y$ at the screen:
$$
O^Ty:=\frac {\partial v^y} {\partial \nu} \qquad \mbox{on \,} \Sigma^T_\sigma .
$$
Notice that the control operator does not change the physical
dimension: $\la W^T f\ra = \la f\ra$, whereas the observation
operator changes the dimension as follows: $\la O^T y\ra =
 \la {\text{time}}\ra \la {\text{length}}\ra^{-1}\la y\ra$.

Let
$$
(y,w)_{\rm int} :=\int_\Omega \frac {\ddd x} {c^2(x)}\, y(x)\,w(x)
$$
and
$$
(f,g)_{\rm ext} :=\int_{\Sigma^T_\sigma} \ddd\Gamma \d t \,\ZZ f(\gamma ,t)
\, g(\gamma ,t)
$$
be the scalar products of functions in $\Omega$ and controls on
the screen. We use the subscript ``int'' in order to emphasize
that we are dealing with functions defined {\it into\/} the
domain, which is not reachable for an external observer, whereas
``ext'' indicates functions on the screen, which the observer can
ope\-ra\-te with. Then (\ref{2.11}) can be written in the form
\begin{equation}\label{2.12}
(W^Tf,y)_{\rm int} =(f,O^Ty)_{\rm ext}
\end{equation}
emphasizing the duality.

A composition of the operators
$$
C^T:=O^T W^T
$$
is an operator mapping controls into observations at the screen.
Thus, one can write
\begin{equation}\label{2.13}
C^Tf=O^T W^Tf=O^Tu^f(\cdot ,T).
\end{equation}
The operator $C^T$ changes the dimension in the same way as $O^T$:
$\la C^T f\ra=\la{\text{time}}\ra\la{\text{length}}\ra^{-1}\la
f\ra$.

\subsection{Jumps in the dual system}

Take a smooth function $y$ and introduce the truncated functions
$$
y_\xi :=\begin{cases}
y & \mbox{in } \Omega^\xi ,\\
0 & \mbox{in } \Omega \setminus \Omega^\xi  \end{cases}
$$
and
$$
y_\xi^\bot :=y-y_\xi =\begin{cases}
0 & \mbox{in } \Omega^\xi ,\\
y & \mbox{in } \Omega \setminus \Omega^\xi . \end{cases}
$$
Note that $y_\xi^\bot$ has a jump at the surface~$\Gamma^\xi$; in
particular, the obvious relation
\begin{equation}\label{2.14}
y_\xi^\bot (x(\gamma ,\xi +0))=y(x(\gamma ,\xi ))
\end{equation}
holds on the part $\sigma^\xi$ lighted by the rays in the tube $B^T$.

Let us insert $y_\xi^\bot$ as the Cauchy data in the second
condition of~(\ref{2.9}). Such a discontinuous perturbation
produces a discontinuous wave carrying a jump at its forward
front. The initial jump of the amplitude (\ref{2.14}) at the point
$x(\gamma,\xi )$ propagates (in inverted time!) along the
ray~$r_\gamma$ and reaches the boundary at the point $x(\gamma
,0)=\NB \gamma$ at the moment $t=T-\NB \xi$ \,(see Fig.~3).

\begin{figure}[htbp]
\centering
\includegraphics {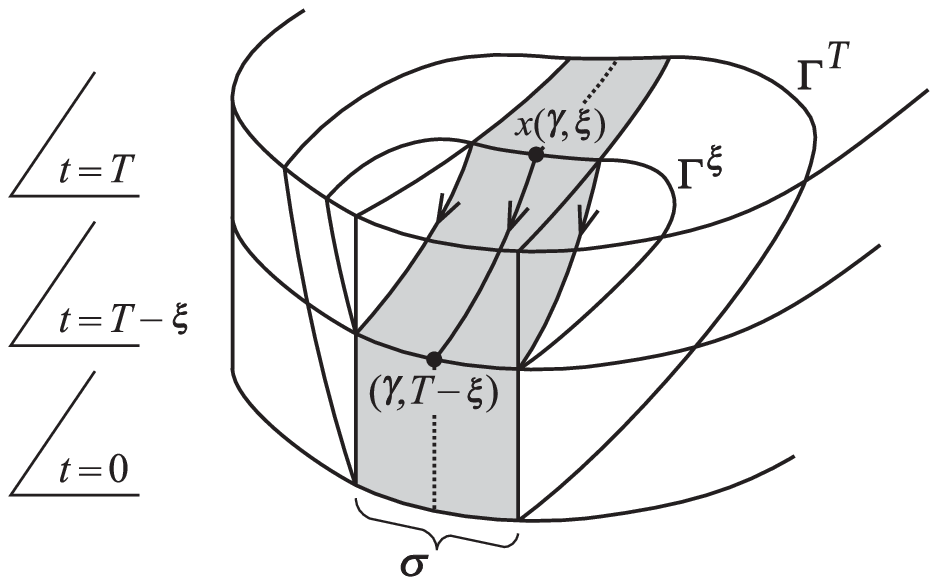}\\[1.5ex]
\small\bf Fig.~3. The jump of $v^{y_\xi^\bot}$
\end{figure}
\noindent Reaching the boundary it produces a jump of the force,
whose amplitude can be calculated by the Geometric Optics formula
\begin{equation}\label{2.15}
\frac {\partial v^{y_\xi^\bot}} {\partial \nu} (\gamma ,T-\xi -0)
=\sqrt { \frac {J(\gamma ,\xi )\,J(\gamma ,0)}
{c(x(\gamma ,\xi ))\,c(x(\gamma ,0)) }}\,
y(x(\gamma ,\xi )).
\end{equation}
This equality can be derived from~(\ref{2.7}) and the duality
relation~(\ref{2.11}) (see, e.g.,~[\ref{BelBlag}]) or in the framework of the
standard ray method~[\ref{BB}]. The obvious duality between (\ref{2.7}) and
(\ref{2.15}) is what physicists call
the {\it reciprocity law}.

Using the observation operator one can rewrite (\ref{2.15}) in the form
\begin{equation}\label{2.16}
[O^T y^\bot_\xi ](\gamma ,T-\xi-0)=
\beta (\gamma ,\xi )\, y(x(\gamma ,\xi ))
\end{equation}
with
$$
\beta (\gamma ,\xi ):=
\sqrt { \frac {J(\gamma ,\xi )\,J(\gamma ,0)}
{c(x(\gamma ,\xi ))\,c(x(\gamma ,0))}}.
$$

\section{Controllability and wave products}
\subsection{Wave shaping}

Choose a function $y$ (of dimension $\la y\ra=\la u^f\ra$) in the
subdomain $\Omega^T$ filled with waves and consider the {\it
boundary control problem\/}: find a control~$f$ providing the
equality
$$
u^f (\cdot ,T)=y \qquad \mbox{in }\, \Omega^T .
$$
In other words, the question is whether one can manage a shape of
waves from the screen. Being of principal character for our
approach, this problem gives the name of the method.

In the general case the answer is the following. For each given $y$ and
arbitrary (small) $\varepsilon >\NB 0$ one can find a control~$f$ such that
$$
\int_{\Omega^T} \frac {\ddd x} {c^2(x)}\, |y(x)-u^f(x,T)|^2 <\varepsilon.
$$
It means that the set of waves is rich enough to approximate any
function (in the mean-square metric). This property of the system
\mbox{(\ref{2.1})--(\ref{2.3})} is called the {\it boundary controllability}.

Controllability is a fact of positive character for inverse
problems: a~very general principle of  system  theory claims that
the richer is the set of states of a dynamical system which an
external observer can create by means of controls, the richer is
information about the system which the observer can extract from
external measurements.

One of the consequences of controllability is the existence of
{\it wave basises}. A~wave basis is a set of waves $u^{f_1}(\cdot
,T)$, $u^{f_2}(\cdot ,T),\dotsc$ satisfying the orthogonality
conditions
$$
\int_{\Omega^T} \frac {\ddd x} {c^2(x)} \,u^{f_i}(x,T)\,
u^{f_j}(x,T)=\delta_{ij}:=
\begin{cases}
1, &i=j\\0,&i\not=j
\end{cases}
$$
such that each $y$ supported in the subdomain~$\Omega^T$ filled
with waves from the screen can be represented in the form of the series
$$
y=\sum_{j=1}^\infty c_j \,u^{f_j} (\cdot ,T)
$$
with the coefficients
$$
c_j =\int_{\Omega^T} \frac {\ddd x} {c^2(x)} \,y(x)\,
u^{f_j}(x,T).
$$
We'll return to wave bases later, after preparing some auxiliary
results.

\subsection{Wave products}
Here we show that scalar products of waves can be expressed via
measurements at the screen.

Take a control $f$ which acts at the screen. Let
$$
f_- (\gamma ,t):=\begin{cases}
f(\gamma, t),   &  0\leqslant t<T, \\
-f(\gamma ,2T-t), & T\leqslant t\leqslant 2T, \end{cases}
$$
be its odd extension and let
$$
\tilde f (\gamma ,t):=\int_0^t f_-(\gamma ,s)\d s \qquad
\mbox{on }\, \Sigma^{2T}.
$$
Take $\tilde f$ as a control in the problem with the doubled final moment:
\begin{align}
\label{3.1}& u_{tt}-c^2\,[\Delta u-qu]=0 \qquad \mbox{in }\, Q^{2T}, \\
\label{3.2}& u|_{t=0}=u_t|_{t=0}=0 \qquad \mbox{in }\, \Omega , \\
\label{3.3}& u={\tilde f} \qquad \mbox{on }\, \Sigma^{2T};
\end{align}
and let $u=u^{\tilde f}(x,t)$ be its solution. We are going to derive the
following important relation
\begin{equation}\label{3.4}
[O^T u^f(\cdot ,T)](\gamma ,t)=
\frac 1 2 \, \Big[
      \frac {\partial u^{\tilde f}} {\partial \nu} (\gamma ,t)-
      \frac {\partial u^{\tilde f}} {\partial \nu} (\gamma ,2T-t)
             \Big]
\end{equation}
for every point $(\gamma ,t)$ of the screen $\Sigma^T_\sigma$.

Consider the auxiliary problem
\begin{align}
\label{3.5}& v_{tt}-c^2\,[\Delta v-qv]=0 \qquad \mbox{in }\, Q^T, \\
\label{3.6}& v|_{t=T}=0,\quad  v_t|_{t=T}=2u_t^{\tilde f}(\cdot ,T) \qquad \mbox{in }\, \Omega ,\\
\label{3.7}& v=0 \qquad \mbox{on }\, \Sigma^T.
\end{align}
Since the coefficients of equation (\ref{3.1}) do not depend on~$t$,
one has
$$
u_t^{\tilde f}(\cdot ,T)=
u^{{\tilde f}_t}(\cdot ,T)=
u^{f_-}(\cdot ,T)=
u^{f}(\cdot ,T);
$$
hence the derivative in (\ref{3.6}) can be calculated as follows:
\begin{equation}\label{3.8}
v_t|_{t=T}=2u_t^{\tilde f}(\cdot ,T)=2u^f(\cdot ,T).
\end{equation}
On the other hand, the solutions of the problems \mbox{(\ref{3.1})--(\ref{3.3})}
and \mbox{(\ref{3.5})--(\ref{3.7})} are connected by the relation
\begin{equation}\label{3.9}
v(x,t)=u^{\tilde f}(x,t)-u^{\tilde f}(x,2T-t) \qquad \mbox{in }\, Q^T
\end{equation}
since the right-hand side of (\ref{3.9}) satisfies all the conditions
\mbox{(\ref{3.5})--(\ref{3.7})}. Differentiating in (\ref{3.9}) one obtains
\begin{equation}\label{3.10}
\frac {\partial v} {\partial \nu} (\gamma ,t)=
\frac {\partial u^{\tilde f}} {\partial \nu} (\gamma ,t)-
\frac {\partial u^{\tilde f}} {\partial \nu}(\gamma ,2T-t)
\qquad \mbox{on }\, \Sigma^T_\sigma .
\end{equation}
Recalling the definition of the observation operator, with regard to
(\ref{3.8}) we get
$$
\{ O^T[2u^f(\cdot ,T)]\} (\gamma ,t)=
\frac {\partial v} {\partial \nu} (\gamma ,t)
\EQF {(\ref{3.10})}
\frac {\partial u^{\tilde f}} {\partial \nu} (\gamma ,t) -
\frac {\partial u^{\tilde f}} {\partial \nu} (\gamma ,2T-t),
$$
which leads to (\ref{3.4}).

In accordance with (\ref{2.13}) the last relation can be written
in the form of representation
\begin{equation}\label{3.11}
(C^Tf) (\gamma ,t)=\frac 1 2 \,
\Big[
\frac {\partial u^{\tilde f}} {\partial \nu} (\gamma ,t)-
\frac {\partial u^{\tilde f}} {\partial \nu}(\gamma ,2T-t)
\Big]
\end{equation}
on the screen.

Choose a pair of controls $f,g$; by virtue of (\ref{3.11}) one has
\begin{gather}
(u^f(\cdot ,T),u^g(\cdot ,T))_{\rm int}
\EQF {(\ref{2.12})}
(O^T u^f (\cdot ,T),g)_{\rm ext}
\EQF {(\ref{2.13})}
(C^T f,g)_{\rm ext} \notag\\
\label{3.12}
{} \EQF {(\ref{3.11})}
\frac 1 2 \int_{\Sigma^T_\sigma} \ddd\Gamma \d t\,
\Big[
\frac {\partial u^{\tilde f}} {\partial \nu} (\gamma ,t)-
\frac {\partial u^{\tilde f}} {\partial \nu}(\gamma ,2T-t)
\Big] \,g(\gamma ,t).
\end{gather}
These relations are interesting and important. An external
observer operating at the boundary can set controls and create
waves but cannot see the waves in themselves. Nevertheless,
(\ref{3.12}) gives him the possibility to calculate the scalar
products of these invisible waves through the data $\partial
u^{\tilde f}\!/ \partial \nu$ and~$g$ on the screen
$\Sigma^T_\sigma$, which are at his disposal.
A.S.Blagovestchenskii was the first who payed attention to this
remarkable fact.

\subsection{Wave basis}

The wave bases introduced at the end of Section~3.1 provide the
main device of the procedure, which makes the waves visible
through the boundary measurements.

Fix $\xi$ provided $0<\xi <T$ and choose a complete system of
controls $g_1, g_2,\dotsc$ supported on the part $\sigma \times
[T-\NB \xi ,T]$ of the screen (so that each $g_j$ is switched on
with delay $T-\NB \xi$ and acts $\xi$~units of time). A
completeness means that controls supported on the same part
$\sigma \times [T-\NB\xi ,T]$ can be expanded as
$f=\sum_{j=1}^\infty \alpha_j g_j$. Notice that in accordance
with~(\ref{2.4}), the corresponding waves $u^{g_j} (\cdot ,T)$ are
supported in the subdomain~$\Omega^\xi$.

The external observer, which possesses the controls $g_j$ and the
measurements $\partial u^{{\tilde g}_j \!} /\partial \nu$, can
modify the system $g_1,g_2,\dotsc$ by the Schmidt process as
follows:
\begin{alignat}{2}
g_1^\prime &:=g_1, \qquad &
f_1 &:=(C^Tg_1^\prime ,g_1^\prime )_{\rm ext}^{-1/ 2}\, g_1^\prime ; \notag\\
g_2^\prime &:=g_2 -(C^Tg_2 ,f_1)_{\rm ext}\, f_1 , \qquad &
f_2 &:=(C^Tg_2^\prime ,g_2^\prime )_{\rm ext}^{-1/2}\, g_2^\prime ; \notag\\
& \hskip 2em \lefteqn{\cdots\cdots\cdots\cdots\cdots\cdots\cdots\cdots
  \cdots\cdots\cdots\cdots\cdots\cdots} \notag\\
g_j^\prime &:=g_j -\sum_{k=1}^{j-1} (C^Tg_j ,f_k)_{\rm ext}\, f_k ,
     \qquad &
f_j &:=(C^Tg_j^\prime ,g_j^\prime )_{\rm ext}^{-1/2}\, g_j^\prime ; \notag\\*
\label{3.13}
& \hskip 2em \lefteqn{\cdots\cdots\cdots\cdots\cdots\cdots\cdots\cdots
  \cdots\cdots\cdots\cdots\cdots\cdots}
\end{alignat}
calculating all the ``ext''-products by (\ref{3.11}) and
(\ref{3.12}). The obtained system $f_1$, $f_2,\dotsc$ is complete
on $\sigma \times [T-\NB \xi ,T]$ and satisfies
$$
(C^T f_i ,f_j)_{\rm ext} =\delta_{ij}.
$$

By virtue of (\ref{3.12}) the corresponding waves
$u^{f_1}(\cdot ,T)$, $u^{f_2}(\cdot ,T),\dotsc$ satisfy
$$
(u^{f_i}(\cdot ,T),u^{f_j}(\cdot ,T))_{\rm int} =
(C^T f_i ,f_j )_{\rm ext}=\delta_{ij} ,
$$
i.e., constitute an orthogonal normalized system in~$\Omega^\xi$.
A~deep fact is that, due to controllability, this system turns out
to be a basis in~$\Omega^\xi$.

Thus, not seeing the waves in themselves, the external observer,
nevertheless, is able to construct a system of controls $f_1$,
$f_2,\dotsc$ producing the wave basis $u^{f_1}(\cdot ,T)$,
$u^{f_2}(\cdot ,T),\dotsc$ in the prescribed
subdomain~$\Omega^\xi$ reachable for waves generated at the
screen.

As a computational problem the orthogonalization by (\ref {3.13})
is equivalent to inversion of the ill-posed Gram matrix
$G_{ij} := (C^T g_i,g_j)_{\rm ext}$ of large size (see~[\ref{BG}]).

The wave basis can be used in the truncation procedure: if $y$ is a function
on~$\Omega$, then its cut-off function $y_\xi$ can be represented as follows:
$$
y_\xi =\sum_{j=1}^\infty (y,u^{f_j}(\cdot ,T))_{\rm int}
\,u^{f_j}(\cdot ,T).
$$
In the important particular case of $y=u^f(\cdot ,T)$ one has
\begin{equation}\label{3.14}
u^f(\cdot ,T)_\xi = \sum_{j=1}^\infty (u^f(\cdot ,T),u^{f_j}(\cdot
,T))_{\rm int} \,u^{f_j}(\cdot ,T) \EQF {(\ref{3.12})}
\sum_{j=1}^\infty (C^Tf,f_j )_{\rm ext} \,u^{f_j}(\cdot ,T)\,,
\end{equation}
so that the external observer can find the coefficients of this
expansion by (\ref{3.12}).

\section{Visualization of waves}
\subsection{Portraits}

Let $y$ be a function in~$\Omega$. A function on~$\Sigma^T_\sigma$
\begin{equation}\label{4.1}
\tilde y(\gamma ,t):=\beta (\gamma ,t)\,y(x(\gamma ,t)),
\end{equation}
($\beta$ is defined in (\ref{2.16})) is called a {\it portrait\/}
(or image) of $y$ on the screen.

Thus, up to the factor $\beta$, the portrait is just the result of
point-wise transferring a function from the tube to the screen.
This is the basic notion of our approach. Note that since $\tilde
y$ is determined by the values $y|_{B^T}$, it would be more
precise to speak about the portrait of the part of $y$ in~$B^T$.

In the rest of the paper, we show that the external observer,
which possesses a complete system of controls $g_1$, $g_2,\dotsc$
on $\Sigma^T_\sigma$ and the corresponding measurements $\partial
u^{\tilde g_1\!}/ \partial \nu$, $\partial u^{\tilde g_2\!}/
\partial \nu,\dotsc$ on $\Sigma^{2T}_\sigma$, is able to reconstruct the
portrait of {\it any\/} wave $u^f(\cdot ,T)$, i.e., can make the
wave to be visible.

\subsection{Amplitude formula}

The procedure for recovering the portraits is based on the Geometric
Optics formulas. Combining the definition (\ref{4.1}) with (\ref{2.16})
one gets the equality
\begin{equation}\label{4.2}
{\tilde y}(\gamma ,\xi )=[O^T y^\bot_\xi ](\gamma ,T-\xi -0),\quad
(\gamma ,\xi )\in \Sigma^T_\sigma,
\end{equation}
which is called the amplitude formula: it represents the portrait
as a collection of amplitudes of the jumps, which are induced by
the truncation of~$y$, pass through the medium, and are detected
at the screen.

Now we are going to insert $y=u^f(\cdot ,T)$ in~(\ref{4.2}). We begin with
a representation of the truncated wave. Applying the observation operator
to (\ref{3.14}) we obtain
$$
O^T u^f(\cdot ,T)_\xi = \sum_{j=1}^\infty (C^T f,f_j)_{\rm ext}\,
O^T u^{f_j}(\cdot ,T) \EQF {(\ref{2.13})} \sum_{j=1}^\infty (C^T
f,f_j)_{\rm ext}\, C^T {f_j}\,,
$$
which yields
$$
O^T u^f(\cdot ,T)_\xi^\bot = O^T [u^f(\cdot ,T) - u^f(\cdot
,T)_\xi ] = C^Tf- \sum_{j=1}^\infty (C^T f,f_j)_{\rm ext}\, C^T
{f_j}\,.
$$
Substituting this in (\ref{4.2}) we arrive at the final formula
\begin{equation}\label{4.3}
[\Widetilde {u^f(\cdot ,T)}] (\gamma ,\xi )= \Big\{ C^Tf -
\sum_{j=1}^\infty (C^T f,f_j)_{\rm ext}\, C^T {f_j} \Big\} (\gamma
,T-\xi -0)\,,
\end{equation}
which is a main device for visualization. It represents the
``slice'' of the portrait corresponding to a fixed~$\xi$ in terms
of the boundary measurements.

\subsection{Visualizing the portraits of waves}

Let us summarize our results. If the external observer can measure
the values $\partial u^{\tilde g\!}/ \partial \nu $ on the screen
$\Sigma^{2T}_\sigma$ for rich enough reserve of controls~$g$, he
is able to recover a portrait of any wave by the following
procedure.

{\bf Step 1} ({\it orthogonalization of controls}).\,\,\, Fix $\xi
:0<\xi <\NB T$ and choose a complete system $g_1$, $g_2,\dotsc$ of
controls supported on the part $\sigma \times [T-\NB \xi ,T]$ of
the screen~$\Sigma^T_\sigma$. Construct the system $f_1$,
$f_2,\dotsc$ (see~(\ref{3.13})) by calculating the
``ext''-products by (\ref{3.12}).

{\bf Step 2} ({\it reconstruction of ``slice''}).\,\,\, Find
$C^Tf_1$, $C^Tf_2,\dotsc$ with the help of~(\ref{3.11}).
Specifying~$f$, recover $[\Widetilde {u^f(\cdot ,T)}](\gamma ,\xi
)$ for all $\gamma \in\NB \sigma $ by means of~(\ref{4.3}).

{\bf Step 3} ({\it reconstruction of the portrait}).\,\,\, Varying
$\xi$ and repeating steps~1 and~2 recover $\Widetilde {u^f(\cdot
,T)}$ on the screen.

\subsection{Recovering the potential}
Assume that the wave speed is constant: $c=\NB 1$, whereas the
(unknown) potential~$q$ is variable. In this case the c-metric is
Euclidean, the rays~$r_\gamma$ are straight lines, and the
factor~$\beta$ entering in the definition of portraits can be
regarded as known. Moreover, we know the relation between the ray
and Cartesian coordinates. Therefore, we can recover the portrait
$\Widetilde {u^f(\cdot ,T)}$ and then find the wave itself at each
point $x(\gamma ,\xi )$ in the tube~$B^T$ by the rule
\begin{equation}\label{4.4}
u^f(x(\gamma ,\xi ),T)=\beta^{-1}(\gamma ,\xi )\,
[\Widetilde {u^f(\cdot ,T)}] (\gamma ,\xi ).
\end{equation}
Finally, specifying $f$ and recovering the waves $u^f(\cdot ,T)$ and
$u^f_{tt} (\cdot ,T)=u^{f_{tt}} (\cdot ,T)$ through their portraits
by~(\ref{4.4}), one can find the potential in the tube $B^T$
from the wave equation
$$
q(x)=
[u^f(x,T)]^{-1} \{ \Delta u^f(x,T) -
u^{f_{tt}}(x,T)\}.
$$

Another way is to take the control $f(\gamma ,t)=\theta  (t-\NB (T-\NB \tau
))$, visualize the corresponding wave $u^f(\cdot ,T)$ in the tube~$B^T$, and
extract the potential from the jump of the wave and its derivatives at the
forward front~$\sigma^\tau$ by means of the well-known Geometric Optics
formulas.

If the potential $q= 0$ and the speed~$c$ is variable and unknown,
we meet another situation: the c-metric is not Euclidean and the
rays~$r_\gamma$ and factor~$\beta$ are {\it unknown}. The
correspondence $(\gamma ,\xi )\leftrightarrow x(\gamma ,\xi )$
between coordinates, which was used for recovering the potential
is also unknown. So, this case requires some additional work.

\subsection{Portraits of harmonic functions}
Let $a=a(x)$ be a harmonic function: $\Delta a(x) =\NB 0$
in~$\Omega$. In the calculations following below, we use the
representation
\begin{equation}\label{4.5}
u^f(x,T)=\int_0^T \ddd t\,\ZZ (T-t)\,u^f_{tt}(x,t)\,,
\end{equation}
which holds due to zero Cauchy data (\ref{2.2}). Also, we assume
that the control~$f$ is smooth and vanishes near $t=\NB 0$ so that
the wave $u^f(x,T)$ is supported in~$\Omega^T$ and vanishes near
its forward front~$\Gamma^T$. Owing to the latter the surface
integrals over the front vanish and we have the equalities
\begin{gather}
(a,u^f(\cdot ,T))_{\rm int} =
\int_\Omega \frac {\ddd x} {c^2(x)}\, a(x)\, u^f(x,T) \notag\\
\EQF {(\ref{4.5})}
\int_\Omega \frac {\ddd x} {c^2(x)} \,a(x)
\int_0^T \ddd t\,\ZZ (T-t) \,u^f_{tt} (x,t) \notag\\
=\int_0^T \ddd t \,\ZZ (T-t)
\int_\Omega {\ddd x}\, \frac {a(x)\, u^f_{tt} (x,t)} {c^2(x)} \notag\\
\EQF {(\ref{2.1}) with $q= 0$}
\int_0^T \ddd t \,\ZZ (T-t)
\int_\Omega \ddd x\, a(x)\, \Delta u^f (x,t) \notag\\
\label{4.6}
= \int_0^T \ddd t \,\ZZ (T-t) \int_\Gamma \ddd\Gamma\,
\Big[ a(\gamma )\frac {\partial u^f} {\partial \nu} (\gamma ,t) -
         \frac {\partial a} {\partial \nu} (\gamma )\,f(\gamma ,t)
\Big].
\end{gather}
Thus, using (\ref{4.6}), the external observer can find the
product of any harmonic function and invisible wave through the
boundary measurements. Note that, in this case, the observer needs
to know $\partial u^f \!/ \partial \nu$ not only on~$\sigma$ but
on a wider part of~$\Gamma$ filled with waves at the moment $t=\NB
T$\footnote{however, see Comments at the end of the paper.}.
Equality (\ref{4.6}) enables one to describe the expansion of the
truncated harmonic function in the wave basis:
\begin{gather}
a_\xi =\sum_{j=1}^\infty \alpha_j u^{f_j}(\cdot ,T) , \notag\\
\label{4.7}
\alpha_j =(a,u^{f_j}(\cdot ,T))_{\rm int}
= \int_{\Sigma^T} \ddd\Gamma \d t \,\ZZ(T-t)\,
\Big[
        a(\gamma )\frac {\partial u^{f_j}} {\partial \nu}(\gamma ,t) -
        \frac {\partial a} {\partial \nu} (\gamma )\,f_j(\gamma ,t)
\Big]
\end{gather}
and then find
\begin{equation}\label{4.8}
O^T a_\xi =\sum_{j=1}^\infty \alpha_j O^T u^{f_j}(\cdot ,T)
\EQF {(\ref {2.13})}
\sum_{j=1}^\infty \alpha_j C^T f_j .
\end{equation}
The latter enables one to visualize the portrait:
\begin{equation}\label{4.9}
\tilde a (\gamma ,\xi )=
\{ O^T [a_T-a_\xi ]\} (\gamma ,T-\xi -0).
\end{equation}
So, the external observer can reconstruct the portrait of {\it any\/}
harmonic function by the scheme:
\begin{itemize}

\item [(i)] for each $\xi :0<\xi \leqslant T$, prepare the system
$f_1,f_2,\dotsc$ on $\sigma \times [T-\NB \xi ,T]$ and then find
$C^Tf_1$, $C^Tf_2,\dotsc$;

\item [(ii)] get the coefficients $\alpha_1,\alpha_2,\dotsc$ by
(\ref{4.7}) and find $O^Ta_\xi$ by~(\ref{4.8});

\item [(iii)]  reconstruct the portrait by (\ref{4.9}).

\end{itemize}

Now we are ready to recover the speed in the tube $B^T$.

\subsection{Recovering the speed}
Introduce a function $1(x)= 1$ in $\Omega$ and the Cartesian
coordinate functions $\pi_1 (x),\allowbreak \pi_2 (x),\allowbreak
\pi_3 (x) : \pi_i(x) =\NB x^i$  for $ x=( x^1,x^2,x^3 )$.  All of
these functions are harmonic; therefore, we can determine their
portraits
$$
\tilde 1(\gamma ,t)=\beta (\gamma ,t), \qquad
\tilde \pi_i (\gamma ,t)=\beta (\gamma ,t)\,\pi_i (x(\gamma ,t))
$$
and then find
$$
x^i(\gamma ,t)=\pi_i(x(\gamma ,t))=
\frac {\tilde \pi_i(\gamma ,t)} {\tilde 1(\gamma ,t)} ,\quad
i=1,2,3.
$$
So, for each point $(\gamma ,t)$ of the screen, we find the point
$x(\gamma,t)$ in the tube. In other words, fixing $\gamma$ and
varying $t$ the external observer can see how the point $x(\gamma
,t)$ moves along the ray~$r_\gamma$ {\it in the interior of\/}
$\Omega$. Then, at last, the observer can find the speed in $B^T$
by the obvious equality
$$
c(x(\gamma ,t))=
\bigg\{ \sum_{i=1}^3
\Big[ \frac {\ddd} {\ddd t} x^i (\gamma ,t) \Big] ^2 \bigg\}^{1/2}.
$$

\section{Comments}
\begin{itemize}
\item There are the versions of the BC-method, which recover
$c|_{B^T}$ via $R^{2T}$ given on $\sigma$ only (i.e., without the
use of measurements outside ~$\sigma$). The first version
[\ref{BelIP}] is a kind of the sample algorithm. A more promoted
variant is proposed in [\ref{BIKS}], section 4.1.

\item As is shown in [\ref{BG}, \ref{BGI}, \ref{BIKS}], the
BC-method can be used as a background of numerical algorithms.
Successful results on numerical testing were obtained by
V.Yu.Gotlib. To our great sorrow, this work had to be suspended
because of his death. Later on, this activity was renewed in
[\ref{BIKS}]; the most promoted results see in [\ref{IBS}].

\item
The problem with the Neumann boundary controls $\partial u/
\partial \nu =\NB f$ on~$\Sigma^T$ may be treated along the same lines.
In this case, the amplitude formula contains one additional
differentiation and looks like that:
$$
\tilde y(\gamma,\xi) = \Big[
\frac {\partial} {\partial t}O^T y^\bot_\xi \Big]
(\gamma, T-\xi-0).
$$
\end{itemize}

\section*{References}

\begin{enumerate}

\item \label {BB} V.~M.~Babich and V.~S.~Buldyrev. {\it
Short-Wavelength Diffraction Theory. Asymptotic Methods.}
Sprin\-ger-Ver\-lag, Berlin, 1991.

\item \label {BelIP} M.~I.~Belishev. Boundary control in
reconstruction of manifolds and metrics (the BC-method). {\it
Inverse Problems\/} (1997) {\bf 13}, No.~5, \mbox {R1--R45}.

\item \label {BelBlag} M.~I.~Belishev and A.~S.~Blagovestchenskii.
{\it Dynamical Inverse Problems of the Wave Propagation Theory.}
St.-Pet\-ers\-burg State University, St.-Pet\-ers\-burg, 1999
(in~Russian).

\item \label {BG} M.~I.~Belishev and V.~Yu.~Gotlib. Dynamical
variant of the BC-method: theory and numerical testing. {\it
J.~Inverse Ill-Posed Problems\/} (1999) {\bf 7}, No.~3, 221--240.

\item \label {BGI} M.~I.~Belishev, V.~Yu.~Gotlib, and
S.~A.~Ivanov. The BC-method in multidimensional spectral inverse
problems: theory and numerical illustrations. {\it ESAIM Control
Optim.\ Calc.\ Var.} (1997) {\bf 2}, 307--327.

\item \label {BIKS} M.I.Belishev, I.B.Ivanov, I.V.Kubyshkin, and
V.S.Semenov. Numerical testing in determination of sound speed
from a part of boundary by the bc-method. {\it J.~Inverse
Ill-Posed Problems\/} (2016) {\bf 23}, No.~5, 221--240.

\item \label {IBS} I.B.Ivanov, M.I.Belishev, V.S.Semenov. The
reconstruction of sound speed in the Marmousi model by the
boundary control method. \\ {http://arxiv.org/abs/1609.07586}.
\end{enumerate}

\end{document}